\documentclass[12pt,letterpaper,english]{article}
\usepackage[T1]{fontenc}
\usepackage[latin9]{inputenc}
\usepackage{latexsym}
\usepackage{mathrsfs}
\usepackage{url}
\usepackage{amsthm}
\usepackage{amsmath}
\usepackage{amssymb}

\makeatletter

\newcommand{\lyxmathsym}[1]{\ifmmode\begingroup\def\b@ld{bold}
  \text{\ifx\math@version\b@ld\bfseries\fi#1}\endgroup\else#1\fi}

\theoremstyle{plain}
\newtheorem{thm}{\protect\theoremname}[section]
  \theoremstyle{definition}
  \newtheorem{defn}[thm]{\protect\definitionname}
  \theoremstyle{plain}
  \newtheorem{lem}[thm]{\protect\lemmaname}
  \theoremstyle{plain}
  \newtheorem{cor}[thm]{\protect\corollaryname}

\date{}

\AtBeginDocument{
  
}

\makeatother

\usepackage{babel}
  \providecommand{\corollaryname}{Corollary}
  \providecommand{\definitionname}{Definition}
  \providecommand{\lemmaname}{Lemma}
\providecommand{\theoremname}{Theorem}

\begin{document}

\title{Constructive Projective Extension of an Incidence Plane}

\author{Mark Mandelkern}
\maketitle
\begin{abstract}
\noindent A standard procedure in classical projective geometry, using
pencils of lines to extend an incidence plane to a projective plane,
is examined from a constructive viewpoint. Brouwerian counterexamples
reveal the limitations of traditional pencils. Generalized definitions
are adopted to construct a projective extension. The main axioms of
projective geometry are verified. The methods used are in accord with
Bishop-type modern constructivism.
\end{abstract}
\noindent {\small MSC: primary 51A45; secondary 03F65 }{\small \par}

\noindent {\small Keywords: Projective extension; Incidence plane;
Constructive mathematics }{\small \par}

\section*{Introduction\label{Section -  Intro} }

In the classical theory of projective geometry, it is a fairly simple
matter to extend an incidence plane to a projective plane; a line
at infinity is added, and pencils of parallel lines become the points
at infinity. A projective plane results; the required projective axioms
are satisfied. In a strictly constructive environment, however, such
an extension presents difficulties due to the indeterminate nature
of arbitrary pencils of lines. 

\emph{Background.} An extension of an incidence plane has been constructed
by Heyting {[}H59{]} and van Dalen {[}D63{]}, using intuitionistic
methods; this work left open the question of the validity of the projective
axiom stating that any two lines have a common point. A recent paper
{[}M11{]} gave a Brouwerian counterexample to demonstrate that in
the Heyting extension the common point axiom is constructively invalid. 

\emph{New extension method.} An extension of an incidence plane will
be constructed using less restrictive definitions, admitting the points
and lines of non-specific character that inevitably emerge in a constructive
setting. The main axioms of projective geometry, including the axiom
that any two lines have a common point, will be verified. 

\emph{Pencils of lines.} Traditionally, a pencil of lines is either
a family of parallel lines, or the family of lines passing through
a given point. Here we use the intrinsic properties of these specific
pencils in adopting a definition which includes pencils with non-specific
properties, such as arise in a constructive setting. Points in the
extended plane will be based on these generalized pencils. 

\emph{Virtual lines. }A central problem in the construction of a projective
extension is the difficulty in determining the nature of an extended
line by means of an object in the original plane. A line in the extended
plane may or may not contain points of the original plane; it is in
general impossible to determine, constructively, which case occurs.
This leads to the concept of a \emph{virtual line,} a set of points
which, if nonvoid, is a line. Virtual lines will be used to construct
points and lines in the extended plane. 

\emph{Bishop-type constructivism. }We follow the constructivist principles
set out by Errett Bishop. Applying these principles when reworking
classical mathematics can have interesting and surprising results.
For the distinctive characteristics of Bishop-type constructivism,
as opposed to intuitionism or recursive function theory, see {[}BR87{]}. 

\emph{Axioms. }There are various approaches to the constructivization
of a classical theory. Bishop's proposal is to find constructive versions
of classical theorems, and to give them constructive proofs. Thus
we adopt no new axioms; we use only constructive versions of the usual
classical axioms for an incidence plane. These axioms are valid on
the real plane $\mathbb{R}^{2},$ taking note of Bishop's thesis that
{}``all mathematics should have numerical meaning''.%
\footnote{\emph{Preface,} pages vii-x, in {[}B67{]}; reprinted as \emph{Prolog,
}pages 1-3, in {[}BB85{]}. %
}

\textit{Logical setting.} This work uses informal intuitionistic logic;
it does not operate within a formal logical system. For the origins
of modern constructivism, and the disengagement of mathematics from
formal logic, see Bishop's Chapter 1, {}``A Constructivist Manifesto'',
in {[}B67{]} or {[}BB85{]}; see also {[}S70{]}, {[}R82{]}, and {[}M85{]}.
Concerning the source of misunderstandings in the mathematical community
as to the methods and philosophy of constructivism, see {[}B65{]}.

\section{Preliminaries\label{> Section 1 - Prelim}}

We assume an incidence plane $\mathscr{G=}(\mathscr{P,L}),$ with
axiom groups \textbf{G} and \textbf{L}, definitions, conventions,
and results from Section 2 of {[}M07{]}. 

\emph{Terminology. }Some terminology and notation used here will be
slightly different from that used in {[}M07{]}. The line through points
$Q$ and $R$ will be denoted $\overline{QR}.$ When the lines $l$
and $m$ are distinct and have a common point, we will say that they
are \emph{intersecting} (rather than {}``nonparallel''); the unique
common point will be denoted by $l\times m.$ Note that the condition
\emph{intersecting} is a primary relation for lines; \emph{parallel}
is its negation. For any set $S,$ the term \emph{nonvoid, }and the
expression $S\neq\varnothing,$\emph{ }are applied in the strict sense:
an element of $S$ has been constructed; it is not sufficient to prove
$\neg(S=\varnothing).$ A distinguished line $l_{0}\in\mathscr{L}$
will be selected. 

\textit{Constructive mathematics.} A characteristic feature of the
constructivist program is meticulous use of the connective {}``or''.
To prove {}``$A\textnormal{ or }B$'' constructively, it is required
that either we prove $A$, or we prove $B$; it is not sufficient
to prove the contrapositive $\neg(\neg A\text{ and }\neg B).$ 

To clarify the methods used here, we give examples of familiar properties
of the real numbers which are constructively \textit{invalid}, and
also properties which are constructively \textit{valid}. The following
classical properties of a real number $c$ are constructively invalid:
{}``Either $c<0$ or $c=0$ or $c>0",$ and {}``If $\neg(c\leq0\ ),$
then $c>0".$ The relation $c>0$ is given a strict constructive definition,
with far-reaching significance. Then, the relation $c\leq0$ is defined
as $\neg(c>0)$. A constructively valid property of the reals is the
\textit{Constructive Dichotomy Lemma}: If $c<d,$ then for any real
number $x$, either $x>c,$ or $x<d.$ This lemma is applied as a
constructive substitute for the constructively invalid \textit{Trichotomy}.
For more details, see {[}B67{]} or {[}BB85{]}. 

\textit{Brouwerian counterexamples.} To determine the specific nonconstructivities
in a classical theory, and thereby to indicate feasible directions
for constructive work, Brouwerian counterexamples are used, in conjunction
with omniscience principles. A \emph{Brouwerian counterexample} is
a proof that a given statement implies an omniscience principle. In
turn, an \emph{omniscience principle} would imply solutions or significant
information for a large number of well-known unsolved problems. This
method was introduced by L. E. J. Brouwer {[}Br08{]} to demonstrate
that use of the \emph{law of excluded middle} inhibits mathematics
from attaining its full significance. Examples will be constructed
on the real plane $\mathbb{R}^{2}.$ For any real number $t$, the
line on $\mathbb{R}^{2}$ with equation $y=t$ will be denoted $l_{t}.$ 

Omniscience principles may be stated in terms of real numbers; we
will have use for the following: 

\medskip{}

\noindent \textbf{\textit{\emph{Limited principle of omniscience (LPO).}}}
\emph{For any real number $c$, either $c=0$ or $c\neq0.$ }

\smallskip{}

\noindent \textbf{\textit{\emph{Weak limited principle of omniscience
(WLPO).}}} \emph{For any real number $c$, either $c=0$ or $\neg(c=0).$}

\smallskip{}

\noindent \textbf{\textit{\emph{Lesser limited principle of omniscience
(LLPO).}}} \emph{For any real number $c,$ either $c\leq0$ or $c\geq0.$ }

\smallskip{}

\noindent \textbf{\textit{\emph{Markov's principle (MP).}}} \emph{For
any real number $c$, if $\neg(c=0)$, then $c\neq0.$ }

\medskip{}

\noindent A statement is considered \emph{constructively invalid}
if it implies an omniscience principle. Following Bishop, we may at
times use the italicized \emph{not} to indicate a constructively invalid
statement. For more information concerning Brouwerian counterexamples,
and other omniscience principles, see {[}B67{]} or {[}BB85{]}, {[}M83{]},
{[}M88{]}, {[}M89{]}, and {[}R02{]}.

\section{Pencils\label{> Section 2 -  Pencils}}

The definition for \emph{pencil of lines} will involve the intrinsic
properties found in pencils of specific type. This will ensure the
inclusion of pencils of unknown type that arise in a constructive
environment. The definition will also admit pencils for which no lines
have been previously constructed, since this situation often occurs
in a constructive setting. We assume the incidence plane $\mathscr{G=}(\mathscr{P,L})$
as indicated in Section \ref{> Section 1 - Prelim}. 
\begin{defn}
\label{Defn Pencil} ~

\noindent $\bullet$ For any point $Q\in\mathscr{P},$ we define 

\noindent 
\[
Q^{*}=\{l\in\mathscr{L}:Q\in l\}
\]

\noindent $\bullet$ For any line $l\in\mathscr{L},$ we define 

\noindent 
\[
l^{*}=\{m\in\mathscr{L}:m\parallel l\}
\]

\noindent $\bullet$ A family of lines $\rho,$ of the form $Q^{*},$
where $Q\in\mathscr{P},$ or $l^{*},$ where $l\in\mathscr{L},$ will
be said to be a \emph{regular} \emph{pencil}. 

\noindent $\bullet$ A family of lines $\alpha$ will be said to be
a \emph{pencil} if it satisfies these two conditions:

(1) $\alpha$ cannot contain fewer than two lines. That is, $\neg(\alpha=\varnothing)$
and $\neg(\alpha\text{ is a singleton)};$ equivalently, $\neg(l,m\in\alpha\text{ implies }l=m).$ 

(2) If $l$ and $m$ are distinct lines in $\alpha$ with $l,m\in\rho,$
where $\rho$ is a regular pencil, then $\alpha\subset\rho$. 

\noindent $\bullet$ A pencil of the form $Q^{*}$ will be said to
be a \emph{point pencil}. 

\noindent $\bullet$ A pencil $\alpha$ with the property that $l\parallel m,$
for any lines $l$ and $m$ in $\alpha,$ will be said to be a \emph{parallel
pencil}. 

\noindent $\bullet$ A pencil $\alpha$ will be said to be \emph{complete}
if the following condition holds:

(2A) If $l$ and $m$ are distinct lines in $\alpha$ with $l,m\in\rho,$
where $\rho$ is a regular pencil, then $\alpha=\rho$. 

\noindent $\bullet$ A pencil $\alpha$ will be said to be \emph{strictly
complete} if the following condition holds:

(2B) If $\alpha\subset\rho,$ where $\rho$ is a regular pencil, then
$\alpha=\rho$.

\noindent $\bullet$ For any pencil $\alpha$, and any line $l$,
we say that $l$ \emph{lies outside} $\alpha,$ written $l\notin\alpha,$
if $l\neq m$ for all lines $m\in\alpha$.

\noindent $\bullet$ Pencils $\alpha$ and $\beta$ are said to be
\emph{distinct}, written $\alpha\neq\beta,$ if there exists a line
$l\in\alpha$ such that $l\notin\beta,$ or there exists a line $l\in\beta$
such that $l\notin\alpha.$
\end{defn}
\noindent \emph{Notes for Definition \ref{Defn Pencil}.} 

(i) \emph{Not} every parallel pencil is regular. For a Brouwerian
counterexample, let $c\in\mathbb{R}$ with $\neg(c=0).$ On $\mathbb{R}^{2},$
set $\alpha=\{l\in\mathscr{L}:c\neq0\text{ and }l\parallel l_{0}\}.$
It is evident that $\alpha$ is a parallel pencil. By hypothesis,
$\alpha=m^{*}$ for some line $m\in\mathscr{L};$ thus $m\in\alpha,$
and $c\neq0.$ Hence MP results. 

For an alternative counterexample, let $c\in\mathbb{R}$ and set  $\beta=\{l\in\mathscr{L}:c=0\mbox{ and }l\parallel l_{0}\}\cup\{l\in\mathscr{L}:c\ne0\text{ and }l\parallel l_{0}\};$
now the hypothesis implies LPO. This example may be easily modified
so that the \emph{Law of Excluded Middle (LEM)} results.

(ii) In the definition of \emph{pencil}, condition (2A) would not
be a suitable substitute for condition (2); this will be indicated
by the Brouwerian counterexample in Note (i) following Theorem \ref{Thm Union class}. 

(iii) Adding condition (2B) to the definition of \emph{pencil} would
complicate the construction of pencils in Theorem \ref{Thm phi(l,m)}
and Theorem \ref{Thm. phi(p,q)}; this will be indicated by Brouwerian
counterexamples in the notes following these theorems.

\title{\noindent }
\begin{lem}
\label{Lm. Pens in Reg 4x} A pencil may be contained in at most one
regular pencil. \end{lem}
\begin{proof}
Let $\alpha$ be any pencil. First assume that $\alpha\subset Q^{*}$
for some point $Q,$ and also $\alpha\subset l^{*}$ for some line
$l.$ Let $m,n\in\alpha$ and suppose that $m\neq n;$ then the lines
$m$ and $n$ intersect at $Q,$ and are also parallel to $l,$ a
contradiction. Thus $m=n,$ and $\alpha$ contains fewer than two
lines, a contradiction. 

Now let $\alpha\subset Q^{*}$ and also $\alpha\subset R^{*},$ for
points $Q$ and $R,$ and suppose that $Q\neq R.$ Then for any line
$l\in\alpha,$ we have $l=\overline{QR},$ so $\alpha$ contains only
one line, a contradiction. Thus $Q=R,$ and $Q^{*}=R^{*}.$ 

Finally, let $\alpha\subset l^{*}$ and also $\alpha\subset m^{*},$
for lines $l$ and $m.$ For any lines $n_{1},n_{2}\in\alpha,$ we
have $n_{1}\parallel l$ and $n_{2}\parallel m.$ Suppose that $l$
and $m$ intersect; it follows from Proposition 2.11 of {[}M07{]}
that $n_{1}$ and $n_{2}$ intersect, so $\alpha$ is contained in
a point pencil, a contradiction. Hence $l\parallel m,$ so $l^{*}=m^{*}.$ \end{proof}
\begin{lem}
Let $Q\in\mathscr{P}$ and let $l,m\in\mathscr{L}.$ 

(a) $l\notin Q^{*}$ if and only if $Q$ lies outside $l$.

(b) $m\notin l^{*}$ if and only if m intersects l.\end{lem}
\begin{proof}
(a) First let $Q\notin l.$ If $m\in Q^{*},$ then $Q\in m,$ so $l\neq m.$
Thus $l\notin Q^{*}.$ Now let $l\notin Q^{*}$. Construct the line
$n$ so that $Q\in n$ and $n\parallel l;$ then $n\in Q^{*},$ so
$l\neq n.$ Thus $l$ and $n$ are parallel and distinct, so $Q\notin l.$

(b) First let $m$ intersect $l,$ and let $n\in l^{*}$. Then $n\parallel l$,
so $m$ intersects $n,$ and $m\neq n$. Thus $m\notin l^{*}$. Now
let $m\notin l^{*}$, choose any point $Q\in m,$ and draw the line
$n$ so that $Q\in n$ and $n\parallel l$. Then $n\in l^{*}$, so
$m\neq n$, and $m$ intersects $n;$ thus $m$ intersects $l.$ \end{proof}
\begin{cor}
\label{regular tite} Let $\rho$ be any regular pencil, and $l$
any line. If $\neg(l\notin\rho),$ then $l\in\rho$. 
\end{cor}
\noindent \emph{Notes for Corollary \ref{regular tite}.} 

(i) For an arbitrary pencil, this property would be constructively
invalid. For a Brouwerian counterexample, let $c\in\mathbb{R}$ with
$\neg(c=0),$ and set $\alpha=\{l_{t}:t\in\mathbb{R}\text{ and }c\ne0\}.$
Then $\neg(l_{0}\notin\alpha),$ but $l_{0}\in\alpha$ would imply
that $c\ne0;$ hence MP results. 

The same pencil $\alpha,$ together with $\beta=l_{0}^{*},$ may be
used to show that the statement {}``\emph{If $\neg(\alpha\neq\beta),$
then $\alpha=\beta$}\,'' is constructively invalid. See, however,
Lemma \ref{Lm-tite-epts-elns}. 

(ii) This property will be extended to a wider class of pencils in
Theorem \ref{Thm. tite}. 
\begin{cor}
\label{Cor Tite alpha rho} Let $\alpha$ be any pencil, and $l$
any line. If $\neg(l\notin\alpha)$ and $\alpha\subset\rho$ for some
regular pencil $\rho,$ then $l\in\rho.$\emph{ }\end{cor}
\begin{thm}
\label{Thm phi(l,m)}There exists a complete pencil containing any
given lines l and m. Define 
\begin{eqnarray*}
\varphi_{0} & = & \{l,m\}\\
\varphi_{1} & = & \{n\in\mathscr{L}:l\times m\in n\}\\
\varphi_{2} & = & \{n\in\mathscr{L}:n\parallel l\parallel m\}\\
\varphi(l,m) & = & \mathbf{\mathbf{\cup}}_{i}\,\varphi_{i}
\end{eqnarray*}

\noindent Then $\varphi(l,m)$ is a complete pencil containing l and
m. \end{thm}
\begin{proof}
Assume that $\varphi(l,m)$ contains fewer than two lines. Suppose
that $l$ intersects $m,$ and set $Q=l\times m;$ then $\varphi_{1}=Q^{*},$
a contradiction. Thus $l\parallel m,$ and $\varphi_{2}=l^{*},$ a
final contradiction. This shows that the assumption is contradictory. 

Let $n_{1}$ and $n_{2}$ be distinct lines in $\varphi(l,m)$ with
$n_{1},n_{2}\in\rho$ for some regular pencil $\rho.$ In the case
$\rho=Q^{*}$ for some point $Q,$ we have $n_{1}\times n_{2}=Q.$
Each of the lines $n_{1}$ and $n_{2}$ lies in one of the sets $\varphi_{i},$
and the required condition $\alpha=Q^{*}$  is easily verified by
considering each of the resulting cases. In the case that $\rho$
is a regular parallel pencil, the verification follows similarly. 

Hence $\varphi(l,m)$  is a complete pencil.
\end{proof}
\noindent \emph{Notes for Theorem \ref{Thm phi(l,m)}. }

(i) It is not assumed that the lines $l$ and $m$ are distinct. 

(ii) The stronger conclusion, {}``$\varphi(l,m)$ is strictly complete'',
would be constructively invalid. For a Brouwerian counterexample,
let $c\in\mathbb{R}$ with $\neg(c=0),$ let $l$ and $m$ be the
lines with equations $y=0$ and $y=cx,$ and consider the pencil $\varphi(l,m).$
Clearly, $\varphi_{0}\subset O^{*},$ where $O$ is the origin. Let
$n\in\varphi_{1}.$ Then $l\neq m,$ and $l\times m=O;$ thus $\varphi_{1}=O^{*},$
and $n\in O^{*}.$ Thus $\varphi_{1}\subset O^{*}.$ If a line is
in $\varphi_{2},$ then $l\parallel m,$ so $c=0,$ a contradiction;
thus $\varphi_{2}=\varnothing$. This shows that $\varphi(l,m)\subset O^{*}.$
By hypothesis, $\varphi(l,m)=O^{*},$ so $m_{0}\in\varphi(l,m),$
where $m_{0}$ is the $y$-axis. It follows that $m_{0}\in\varphi_{1},$
and thus $c\neq0.$ Hence MP results. 

(iii) When $l\neq m,$ the Heyting {}``p.point'' is defined in {[}H59{]}
by $\mathfrak{P}(l,m)=\{n\in\mathscr{L}:n\cap l=l\cap m\text{ or }n\cap m=l\cap m\}.$
Clearly, $\varphi(l,m)\subset\mathfrak{P}(l,m),$ and $\neg(\varphi(l,m)\neq\mathfrak{P}(l,m))$.
However, equality here would be constructively invalid. For a Brouwerian
counterexample, let $c\in\mathbb{R}$ and let $l$ and $m$ be the
lines with equations $y=0$ and $y=1-cx.$ Consider the line with
equation $y=2-2cx;$ LPO results. 

Many other results in the present paper also have analogues in {[}H59{]}.
However, the basic definitions differ in quite fundamental ways, so
a detailed comparison of the two approaches is not feasible. 

(iv) This method will be extended in Theorem \ref{Thm. phi(p,q)}. 
\begin{lem}
~

(a) A complete pencil $\alpha$ is a point pencil if and only if $\alpha\cap l^{*}\neq\varnothing$
for every regular parallel pencil $l^{*}.$ 

(b) Point pencils $Q^{*}$ and $R^{*},$ if distinct, have the line
$\overline{QR}$ in common. \end{lem}
\begin{proof}
(a) The necessity is the parallel postulate. For the converse, choose
any intersecting lines $m_{1}$ and $m_{2},$ and select lines $n_{1}\in\alpha\cap m_{1}^{*}$
and $n_{2}\in\alpha\cap m_{2}^{*}.$ It follows that $n_{1}$ intersects
$n_{2},$ so $\alpha=Q^{*},$ where $Q=n_{1}\times n_{2}.$ 

(b) We may select a line $l\in R^{*}$ such that $l\notin Q^{*}.$
Then $Q\notin l$ and $R\in l,$ so $Q\ne R.$ Thus $\overline{QR}$
is a line common to both pencils. \end{proof}
\begin{lem}
\label{Lm not regular}For any given complete pencil $\alpha,$ the
statement {}``$\neg(\alpha=\rho)$ for every regular pencil $\rho$'',
is contradictory. However, the statement {}``Every strictly complete
pencil is either a point pencil or a regular parallel pencil'' is
constructively invalid. \end{lem}
\begin{proof}
Assume the first statement. Let $l,m\in\alpha,$ and suppose that
$l\neq m.$ If $l$ intersects $m,$ with $Q=l\times m,$ then $\alpha=Q^{*},$
a contradiction. Thus $l\parallel m,$ so $\alpha=l^{*},$ also a
contradiction. Hence $l=m.$ This shows that $\alpha$ contains fewer
than two lines, a final contradiction. 

For a Brouwerian counterexample to the second statement, let $c\in\mathbb{R},$
let $l$ and $m$ be the lines with equations $y=0$ and $y=1-cx,$
and set $\beta=\varphi(l,m)$ using the construction of Theorem \ref{Thm phi(l,m)}.
The pencil $\beta$ is strictly complete; apply the hypothesis to
$\beta,$ and LPO results. \end{proof}
\begin{thm}
Let $Q^{*}$ be any point pencil, and $\beta$ any complete pencil.
Then $\neg(Q^{*}\cap\beta=\varnothing).$ However, the stronger conclusion,
{}``\textup{$Q^{*}\cap\beta\ne\varnothing$''}, is constructively
invalid.\end{thm}
\begin{proof}
Assume\emph{ }that $Q^{*}\cap\beta=\varnothing.$ Let $l,m\in\beta$
and suppose that $l\neq m.$ Suppose further that $l$ intersects
$m,$ and set $R=l\times m;$ then $\beta=R^{*}.$ Now suppose even
further that $R\ne Q;$ then $\overline{QR}\in Q^{*}\cap\beta,$ a
contradiction. Thus $R=Q,$ and $\beta=Q^{*},$ a contradiction. Thus
$l\parallel m$ and $\beta=l^{*}.$ Since the parallel postulate provides
a line in $l^{*}$ passing through $Q,$ this is a contradiction.
Thus $l=m.$ This shows that $\beta$ contains fewer than two lines,
a final contradiction. Hence $\neg(Q^{*}\cap\beta=\varnothing).$ 

For a Brouwerian counterexample to the stronger conclusion, let $c\in\mathbb{R}$
and take $Q$ at the origin in $\mathbb{R}^{2}.$ Let $l$ and $m$
be the lines with equations $y=1$ and $y=1+c-cx$, and let $\beta=\varphi(l,m)$
be the complete pencil constructed using Theorem \ref{Thm phi(l,m)}.
Set $R=(1,1).$ Note that if $c=0,$ then $\beta=l^{*},$ so the $x$-axis
$l_{0}$ is the unique line common to $Q^{*}$ and $\beta.$ On the
other hand, if $c\ne0,$ then $\beta=R^{*},$ and $\overline{QR}$
is the unique common line. By hypothesis, $Q^{*}$ and $\beta$ have
a common line $n;$ either $n\neq l_{0},$ or $n\neq\overline{QR}.$
It follows that either $\neg(c=0),$ or $c=0.$ Hence WLPO results. 
\end{proof}
\noindent Applications of this theorem will be by way of the following: 
\begin{cor}
\label{Cor parpen} If $\alpha$ is a complete pencil with the property
that $\alpha\cap\beta=\varnothing$ for some complete pencil $\beta,$
then $\alpha$ is a parallel pencil. 
\end{cor}
\noindent The definition of \emph{pencil} includes multiple pencils
which will represent the same point in the projective extension. Thus
we adopt the following: 
\begin{defn}
~~ 

\noindent $\bullet$ Pencils $\alpha$ and $\beta$ will be said to
be \emph{equivalent,} written $\alpha\approx\beta,$ if for any regular
pencil $\rho,$ the following condition is satisfied: $\alpha\subset\rho$
if and only if $\beta\subset\rho$. 

\noindent $\bullet$ The pencils $\alpha$ and $\beta$ are said to
be \emph{non-equivalent} if $\neg(\alpha\approx\beta).$

\noindent $\bullet$ With regard to the resulting equivalence relation
$\approx,$ the equivalence class containing the pencil $\alpha$
will be denoted $\overline{\alpha}$. \end{defn}
\begin{thm}
\label{Thm Pens } Let $\alpha$ and $\beta$ be any pencils.

(a) If $\alpha$ and $\beta$ have two distinct lines in common, then
$\alpha\approx\beta.$ Non-equivalent pencils have at most one line
in common. 

(b) If $\alpha$ is a parallel pencil, and $\alpha\approx\beta,$
then $\beta$ is also a parallel pencil. 

(c) If $\neg(\alpha\neq\beta),$ then $\alpha\approx\beta.$ Conversely,
if $\alpha\approx\beta,$ with $\alpha$ complete and $\beta$ strictly
complete, then $\neg(\alpha\ne\beta).$ \end{thm}
\begin{proof}
(a) Let $l$ and $m$ be lines common to $\alpha$ and $\beta,$ with
$l\neq m.$ If $\alpha\subset\rho$ for some regular pencil $\rho,$
then $l,m\in\rho,$ so $\beta\subset\rho.$ Thus $\alpha\approx\beta.$ 

(b) Let $l,m\in\beta,$ and suppose that $l$ and $m$ intersect at
some point $Q.$ Then $\beta\subset Q^{*},$ so also $\alpha\subset Q^{*}.$
For any lines $u,v\in\alpha,$ we have $u\parallel v.$ If $u\ne v,$
then $\alpha\subset u^{*},$ contradicting Lemma \ref{Lm. Pens in Reg 4x}.
Thus $u=v,$ and $\alpha$ contains fewer than two lines, a contradiction;
hence $l\parallel m.$ 

(c) Let $\neg(\alpha\neq\beta).$ Let $\alpha\subset\rho,$ let $l\in\beta,$
and suppose that $l\notin\rho.$ Then $l\notin\alpha,$ so $\alpha\neq\beta,$
a contradiction. Thus $\neg(l\notin\rho),$ so by Corollary \ref{regular tite}
we have $l\in\rho.$ Thus $\beta\subset\rho.$ This shows that $\alpha\approx\beta.$ 

Now let $\alpha\approx\beta,$ with $\alpha$ complete and $\beta$
strictly complete, and assume that $\alpha\ne\beta.$ Suppose that
$\alpha=\rho$ for some regular pencil $\rho;$ then $\alpha\subset\rho,$
so $\beta\subset\rho,$ and $\beta=\rho,$ a contradiction. Thus $\neg(\alpha=\rho)$
for all regular pencils $\rho.$ By Lemma \ref{Lm not regular}, this
is a contradiction. Hence $\neg(\alpha\ne\beta).$ \end{proof}
\begin{lem}
\label{Reg vs Equiv} ~

(a) Let $\alpha$ be any pencil, and $\rho$ a regular pencil. Then
$\alpha\approx\rho$ if and only if $\alpha\subset\rho.$ 

(b) There exists at most one regular pencil in any given equivalence
class of pencils. \end{lem}
\begin{proof}
In (a), let $\alpha\subset\rho;$ then it follows from Lemma \ref{Lm. Pens in Reg 4x}
that $\alpha\approx\rho.$ The converse is immediate. Now (b) follows
from (a). \end{proof}
\begin{thm}
\label{Thm Par Pens } Let $\alpha$ and $\beta$ be parallel pencils.

(a) If $\alpha$ and $\beta$ have a line in common, then $\alpha\approx\beta.$
Non-equivalent parallel pencils have no common lines. 

(b) If $\alpha$ and $\beta$ have a line in common, and are strictly
complete, then $\alpha=\beta.$ \end{thm}
\begin{proof}
Select a common line $l;$ then $\alpha\subset l^{*}$ and $\beta\subset l^{*}.$
In (a), it follows from Lemma \ref{Reg vs Equiv} that $\alpha\approx\beta.$
In (b), we have $\alpha=l^{*}$ and $\beta=l^{*};$ thus $\alpha=\beta.$\end{proof}
\begin{defn}
For any pencil $\alpha,$ we define
\[
\alpha'=\cup\{\beta:\beta\in\overline{\alpha}\}
\]
After verification in the next theorem, the pencil $\alpha'$ will
be called the \emph{full} pencil\emph{ }in the equivalence class $\overline{\alpha}.$
A pencil $\alpha$ will be said to be \emph{prime} if $\alpha=\alpha'.$ \end{defn}
\begin{thm}
\label{Thm Union class} The union $\alpha'$ of all pencils in an
equivalence class $\overline{\alpha}$ is a strictly complete pencil,
equivalent to the pencils in the class. 

The pencil $\alpha'$ is unique in this limited sense: If $\beta$
is any complete pencil in the equivalence class $\overline{\alpha},$
then $\neg(\beta\ne\alpha').$ However, the statement {}``If $\beta$
is any strictly complete pencil in $\overline{\alpha},$ then $\beta=\alpha'$''
is constructively invalid. 

If $\overline{\alpha}$ contains a regular pencil $\rho,$ then $\alpha'=\rho.$\end{thm}
\begin{proof}
(i) To show that $\alpha'$ is a pencil, let $l$ and $m$ be distinct
lines in $\alpha'$ with $l,m\in\rho$ for some regular pencil $\rho,$
and select pencils $\beta_{1},\beta_{2}\in\overline{\alpha}$ with
$l\in\beta_{1}$ and $m\in\beta_{2}.$ 

Set $\beta=\beta_{1}\cup\beta_{2}.$ To show that $\beta$ is a pencil,
let $n_{1}$ and $n_{2}$ be distinct lines in $\beta$ with $n_{1},n_{2}\in\sigma$
for some regular pencil $\sigma.$ If $n_{1}$ and $n_{2}$ both belong
to $\beta_{1}$, then $\beta_{1}\subset\sigma,$ so also $\beta_{2}\subset\sigma,$
and $\beta\subset\sigma.$ Similarly for $\beta_{2}.$ Thus we may
assume that $n_{1}\in\beta_{1}$ and $n_{2}\in\beta_{2}.$ Let $n\in\beta;$
we may assume that $n\in\beta_{1}.$ Suppose that $n\notin\sigma;$
it follows that $n\ne n_{1},$ with $n,n_{1}\in\beta_{1}.$ Then $\beta_{1}\subset\sigma,$
so also $\beta_{2}\subset\sigma,$ and $\beta\subset\sigma,$ so $n\in\sigma,$
a contradiction. Thus $n\in\sigma.$ This shows that $\beta\subset\sigma.$
Hence $\beta$ is a pencil. Clearly, $\beta\approx\beta_{1},$ so
$\beta\in\overline{\alpha}.$ 

Since $l$ and $m$ are distinct lines in $\beta$ with $l,m\in\rho,$
it follows that $\beta\subset\rho.$ For any pencil $\gamma\in\overline{\alpha},$
we have $\gamma\approx\beta,$ so $\gamma\subset\rho;$ thus $\alpha'\subset\rho.$
Hence $\alpha'$ is a pencil. 

(ii) Let $\rho$ be a regular pencil with $\alpha\subset\rho.$ Then
$\gamma\subset\rho$ for all pencils $\gamma\in\overline{\alpha},$
so $\alpha'\subset\rho.$ Hence $\alpha'\approx\alpha.$ 

(iii) If $\alpha'\subset\rho$ for some regular pencil $\rho,$ then
it follows from Lemma \ref{Reg vs Equiv} that $\alpha'\approx\rho,$
so $\rho\in\overline{\alpha},$ and $\rho\subset\alpha';$ thus $\alpha'=\rho.$
This shows that $\alpha'$ is a strictly complete pencil. 

(iv) If $\beta$ is a complete pencil in the equivalence class $\overline{\alpha},$
then it follows from Theorem \ref{Thm Pens } that $\neg(\beta\ne\alpha').$ 

(v) For a Brouwerian counterexample involving the statement in quotes,
let $c\in\mathbb{R},$ let $m$ be the line on $\mathbb{R}^{2}$ with
equation $y=1-cx,$ and set $\alpha=\varphi(l_{0},m).$ If $\alpha\subset n^{*}$
for some line $n,$ then $c=0$ and $\alpha=l_{0}^{*}=n^{*}.$ If
$\alpha\subset R^{*}$ for some point $R,$ with $R=(d,e),$ then
$e=0$ and $1-cd=0,$ so $c\neq0$ and $\alpha=(1/c,0)^{*}=R^{*}.$ 

Now consider the family of lines $\beta=\{l\in\mathscr{L}:c=0\text{ and }l\parallel l_{0}\}\cup\{l\in\mathscr{L}:c\ne0\text{ and }(1/c,0)\in l\}.$
It is clear that $\beta$ is a pencil. Let $\beta\subset n^{*}$ for
some line $n;$ then $c=0$ and $\beta=l_{0}^{*}=n^{*}.$ Now let
$\beta\subset R^{*}$ for some point $R=(d,e).$ Suppose that $|c|<1/(|d|+1).$
Suppose further that $c\ne0;$ then $\beta=(1/c,0)^{*}.$ It follows
that $R=(1/c,0),$ so $d=1/c$ and $|d|>|d|+1,$ an absurdity. Thus
$c=0,$ and $\beta=l_{0}^{*},$ a contradiction. This shows that $|c|\geq1/(|d|+1),$
so $c\ne0,$ and $\beta=(1/c,0)^{*}=R^{*}.$ 

Thus $\beta$ is strictly complete, and $\beta\approx\alpha.$ By
hypothesis, $\beta=\alpha'.$ Since $l_{0}\in\alpha\subset\alpha',$
it follows that $l_{0}\in\beta,$ and LPO results. 

(vi) If $\rho\in\overline{\alpha}$ for some regular pencil $\rho,$
then $\alpha'\approx\rho,$ so by Lemma \ref{Reg vs Equiv} we have
$\alpha'\subset\rho.$ Thus $\alpha'=\rho.$ 
\end{proof}
\noindent \emph{Notes for Theorem \ref{Thm Union class}.} 

(i) The proof shows that \emph{The union of any two equivalent pencils
is also a pencil, equivalent to the given pencils.} However, the union
of two complete equivalent pencils need \emph{not} be complete. For
a Brouwerian counterexample, let $c\in\mathbb{R}$ with $\neg(c=0).$
On $\mathbb{R}^{2},$ set $\alpha=\{l_{0}\}\cup\{l_{t}:t\in\mathbb{R}\text{ and }c\ne0\}$
and $\beta=\{l_{1}\}\cup\{l_{t}:t\in\mathbb{R}\text{ and }c\neq0\};$
these pencils are complete, and equivalent. The hypothesis applied
to $\alpha\cup\beta$ results in MP. This is one of the facts concerning
complete pencils that necessitates the more inclusive definition of
\emph{pencil} adopted in Definition \ref{Defn Pencil}.

(ii) The union of two strictly complete equivalent pencils is strictly
complete. 

(iii) The counterexample in the proof of the theorem involves an equivalence
class containing \emph{no} regular pencil. On the other hand, it follows
from Lemma \ref{Lm not regular} that the statement {}``If $\rho$
is any regular pencil, then $\neg(\alpha'=\rho)$'' is contradictory. 
\begin{thm}
\label{Thm. tite} Let $\alpha$ be a prime pencil, and $l$ any line.
If $\neg(l\notin\alpha),$ then $l\in\alpha.$ \end{thm}
\begin{proof}
Set $\beta=\alpha\cup\{l\}.$ To show that $\beta$ is a pencil, let
$m_{1}$ and $m_{2}$ be distinct lines in $\beta$ with $m_{1},m_{2}\in\rho$
for some regular pencil $\rho.$ Either both lines are in $\alpha,$
or one line is $l$. 

In the first case, since $\alpha$ is a pencil we have $\alpha\subset\rho.$
It follows from Corollary \ref{Cor Tite alpha rho} that $l\in\rho;$
thus $\beta\subset\rho.$ 

In the second case, we may say that $m_{1}\in\alpha$ and $m_{2}=l.$
Let $n$ be any line in $\alpha,$ and suppose that $n\notin\rho.$
If $\alpha$ were a regular pencil, then by Corollary \ref{regular tite}
we would have $m_{2}=l\in\alpha,$ and it would follow that $\alpha\subset\rho,$
so $n\in\rho,$ a contradiction. Thus $\neg(\alpha\text{ is regular}),$
and by Lemma \ref{Lm not regular} this is contradictory. Thus $n\in\rho.$
This shows that $\alpha\subset\rho,$ and thus $\beta\subset\rho.$ 

Hence $\beta$ is a pencil. To show that $\beta$ is equivalent to
$\alpha,$ let $\alpha\subset\rho$ for some regular pencil $\rho.$
It follows from Corollary \ref{Cor Tite alpha rho} that $l\in\rho,$
and this shows that $\beta\subset\rho.$ Thus $\beta\approx\alpha.$
Since $\alpha$ is a prime pencil, it follows that $\beta\subset\alpha,$
and therefore $l\in\alpha.$ 
\end{proof}

\section{Virtual lines\label{> Section 3 - V-lines} }

Consider the following classical situation: If a line $L$ in the
extended plane contains a proper point, then the set $S,$ of all
proper points on $L,$ is a line in the original plane. However, if
$L$ is the line at infinity, then the set $S$ is void. Constructively,
we will not know in general which is the case. Thus we adopt the following: 
\begin{defn}
~

\noindent $\bullet$ A set $p$ of points in $\mathscr{P}$ will be
said to be a \emph{virtual line, }or\emph{ v-line,} if it satisfies
this condition: If $p$ is nonvoid, then $p$ is a line. 

\noindent $\bullet$ The family of all v-lines will be denoted $\mathscr{V}.$ 

\noindent $\bullet$ The v-lines $p$ and $q$ will be said to be
\emph{distinct}, written $p\neq q,$ if they satisfy these two conditions: 

(1) $\neg(p=q)$ 

(2) If $p,q\in\mathscr{L},$ then $p\neq q$ in the sense of distinct
lines in $\mathscr{L}.$

\noindent $\bullet$ The expression $p\neq\varnothing$ will imply
that the v-line $p$ is nonvoid, and thus $p$ is a line. 

\noindent $\bullet$ When we write $p\times q=Q,$ or $p\parallel q,$
this will imply that $p$ and $q$ are lines. 
\end{defn}
\noindent The notion of \emph{v-line} also arises in connection with
pencils. Theorem \ref{Thm Pens } shows that non-equivalent pencils
have at most one line in common, and Theorem \ref{Thm Par Pens }
shows that non-equivalent parallel pencils have no common lines. Thus
the family of lines common to two non-equivalent pencils may consist
of a single line, or it may be void; constructively, it is in general
unknown which alternative holds. 
\begin{defn}
 For any non-equivalent pencils $\alpha$ and $\beta,$ we define
\[
\alpha\sqcap\beta=\{Q\in\mathscr{P}:Q\in l\in\alpha\cap\beta\text{ for some }l\in\mathscr{L}\}
\]

\noindent The set of points $\alpha\sqcap\beta$ will be called the
\emph{core }of the pair\emph{ }$\alpha,\beta.$ \end{defn}
\begin{lem}
\label{Lm. core is quasi-line} The core \textup{$\alpha\sqcap\beta,$}
where $\alpha$ and $\beta$ are any non-equivalent pencils, \textup{\emph{is
a v-line. }}\end{lem}
\begin{proof}
Set $p=\alpha\sqcap\beta,$ and let $p\ne\varnothing.$ Construct
a point $R\in p$ and a line $m$ such that $R\in m\in\alpha\cap\beta;$
clearly $m\subset p.$ For any point $Q\in p,$ we have $Q\in l\in\alpha\cap\beta$
for some line $l;$ it follows from Theorem \ref{Thm Pens } that
$l=m,$ so $Q\in m.$ This shows that $p\subset m.$ Hence $p=m,$
so $p$ is a line. \end{proof}
\begin{thm}
\label{Thm. phi(p,q)} Let p and q be any v-lines. Define 
\begin{eqnarray*}
\varphi_{0} & = & \{p,q\}\cap\mathscr{L}\\
\varphi_{1} & = & \{l\in Q^{*}:p\times q=Q\}\\
\varphi_{2} & = & \{l\in p^{*}:p\parallel q\}\\
\varphi_{3} & = & \{l\in p^{*}:p\neq\varnothing,q=\varnothing\}\cup\{l\in q^{*}:q\neq\varnothing,p=\varnothing\}\\
\varphi_{4} & = & \{l\in l_{0}^{*}:p=q=\varnothing\}\\
\varphi(p,q) & = & \cup_{i}\,\varphi_{i}
\end{eqnarray*}

\noindent Then $\varphi(p,q)$ is a complete pencil.
\end{thm}
\noindent To aid the proof, and for later use, we have first: 
\begin{lem}
\label{Lemma for phi} In Theorem \ref{Thm. phi(p,q)}, 

(a) If $p\times q=Q,$ then $\varphi(p,q)=Q^{*}.$

(b) If $p\parallel q,$ then $\varphi(p,q)=p^{*}.$

(c) If there exists a line l in $\varphi_{2}\cup\varphi_{3}\cup\varphi_{4},$
then $\varphi(p,q)=l^{*}.$

(d) If $p=\varnothing$ or $q=\varnothing,$ then $\varphi(p,q)$
is a parallel pencil. \end{lem}
\begin{proof}
The first three properties are immediate. For (d), let $q=\varnothing,$
and let $l,m\in\varphi(p,q);$ then $l\in\varphi_{0}\cup\varphi_{3}\cup\varphi_{4}.$
If $l\in\varphi_{0},$ then $l=p,$ so $p$ is nonvoid, and $\varphi_{3}=p^{*}.$
Since $m\in\varphi_{0}$ or $m\in\varphi_{3},$ it follows that $m\parallel l.$
If $l$ lies in $\varphi_{3}$ or $\varphi_{4},$ then (c) applies. 
\end{proof}
\noindent \emph{Proof of Theorem \ref{Thm. phi(p,q)}.} Assume that
$\lyxmathsym{\textquotedblleft}l,m\in\varphi(p,q)$ implies $l=m\lyxmathsym{\textquotedblright}.$
Suppose that one of the v-lines, say $p,$ is nonvoid. Now suppose
that $q\neq\varnothing,$ and then suppose further that $p$ intersects
$q,$ with $Q=p\times q.$ Then $\varphi(p,q)=Q^{*},$ a contradiction.
Thus $p\parallel q,$ so $\varphi(p,q)=p^{*},$ a contradiction. Thus
$q=\varnothing,$ and $\varphi(p,q)=p^{*},$ a contradiction. Thus
both v-lines are void, so $\varphi(p,q)=l_{0}^{*},$ contradicting
the assumption. This shows that $\varphi(p,q)$ cannot contain fewer
than two lines. 

Let $l,m\in\varphi(p,q)$ with $l\neq m.$ First let $l$ intersect
$m,$ with $Q=l\times m;$ then $l,m\in\varphi_{0}\cup\varphi_{1}.$
In the first case, where $l,m\in\varphi_{0},$ it follows that $p\times q=Q,$
so $\varphi(p,q)=Q^{*}.$ In the other three cases, $\varphi_{1}\neq\varnothing,$
and again $\varphi(p,q)=Q^{*}.$ Now let $l\parallel m.$ Using part
(c) of the lemma we may assume that $l,m\in\varphi_{0}\cup\varphi_{1}.$
Since $\varphi_{1}=\varnothing,$ we may say that $l=p$ and $m=q.$
Now $\varphi(p,q)=\varphi_{2}=p^{*}=l^{*}.$ 

Hence $\varphi(p,q)$ is a complete pencil.\hspace*{\fill} $\square$

\bigskip{}

\noindent \emph{Notes for Theorem \ref{Thm. phi(p,q)}.} ~

(i) It is not assumed that the v-lines $p$ and $q$ are distinct. 

(ii) The stronger conclusion, \emph{{}``$\varphi(p,q)$ is} \emph{strictly
complete'',} would be constructively invalid. For a Brouwerian counterexample,
let $c\in\mathbb{R}$ with $\neg(c=0),$ set $p=\{(t,0):t\in\mathbb{R}\text{ and }c\neq0\},$
and set $q=\{(0,t):t\in\mathbb{R}\text{ and }c\neq0\}.$ Then $\varphi(p,q)\subset(0,0)^{*},$
but $\varphi(p,q)=(0,0)^{*}$ would imply MP. 

\bigskip{}

\noindent The definition of \emph{v-line} includes multiple v-lines
which will form the basis for the same line in the projective extension.
Thus we adopt the following: 
\begin{defn}
~ 

\noindent $\bullet$ The v-lines \emph{p} and \emph{q} will be said
to be \emph{equivalent,} written $p\approx q,$ if $\neg(p\neq q).$ 

\noindent $\bullet$ After verification in Theorem \ref{Thm - Eqclasses},
the equivalence class containing the v-line $p$ will be denoted $\overline{p}.$ 

\noindent $\bullet$ For any v-line $p\in\mathscr{V},$ we define

\noindent 
\[
p'=\cup\{q:q\in\overline{p}\}
\]

\noindent $\bullet$ The v-line $p$ will be said to be \emph{prime}
if $p=p'.$ 

\noindent $\bullet$ The family of all prime v-lines will be denoted
$\mathscr{V'}.$ \end{defn}
\begin{lem}
\label{Lm Equiv condition}For any v-lines p and q, the following
conditions are equivalent:

(a) \textup{$p\approx q.$ }

(b) For any line $l\in\mathscr{L},$ $p\subset l$ if and only if
$q\subset l.$\end{lem}
\begin{proof}
Let $p\approx q,$ and let $p\subset l$ for some line $l.$ Let $Q\in q,$
and suppose that $Q\notin l;$ then $\neg(Q\in p),$ so $\neg(p=q).$
If $p\in\mathscr{L},$ then $p=l,$ so $Q\notin p;$ thus $p$ and
$q$ are distinct as lines. This shows that $p$ and $q$ are distinct
as v-lines, a contradiction; thus $Q\in l.$ Hence $q\subset l.$ 

Let condition (b) hold, and assume that $p\neq q.$ Under these conditions,
it is clear that the two v-lines cannot both be nonvoid. Suppose that
one of the v-lines, say $p,$ is nonvoid; it follows that $q$ is
void. Select any two distinct lines; since $q$ is contained in both
lines, condition (b) implies that the line $p$ is contained in both
lines, an absurdity. Thus both v-lines are void. Now $p=q,$ a contradiction.
This shows that $\neg(p\neq q),$ so $p\approx q.$ \end{proof}
\begin{thm}
\label{Thm - Eqclasses} The relation $\approx$ on the family $\mathscr{V}$
of v-lines is an equivalence relation. The union $p'$ of all v-lines
in an equivalence class $\overline{p}$ is a v-line, equivalent to
the v-lines in the class. \end{thm}
\begin{proof}
It follows from the lemma that $\approx$ is an equivalence relation.
Let $p'\neq\varnothing,$ and select a v-line $r\in\overline{p}$
with $r\neq\varnothing;$ thus $r$ is a line. For any v-line $q\in\overline{p},$
it follows that $q\subset r.$ Thus $p'\subset r,$ so $p'=r.$ Hence
$p'$ is a v-line. If $p\subset l$ for some line $l,$ then $q\subset l$
for all $q\in\overline{p};$ thus $p'\subset l.$ This shows that
$p'\approx p.$ \end{proof}
\begin{lem}
\label{v-line class}Let $\overline{p}$ be an equivalence class of
v-lines. 

(a) If q and r are nonvoid v-lines in $\overline{p},$ then $q=r.$ 

(b) If there exists a v-line q in $\overline{p}$ with $q=\varnothing,$
then $r=\varnothing$ for all \textup{$r\in\overline{p}.$ }\end{lem}
\begin{proof}
(a) Since $q$ and $r$ are both lines, and $\neg(q\neq r),$ it follows
from Proposition 2.16 in {[}M07{]} that $q=r.$ 

(b) Let $r\in\overline{p},$ and suppose that $r\neq\varnothing;$
thus $r$ is a line. Select any two distinct lines; since $q$ is
contained in both lines, it follows from Lemma \ref{Lm Equiv condition}
that $r$ is contained in both lines, an absurdity. Thus $r=\varnothing.$ 
\end{proof}

\section{Extension points and extension lines\label{> Section 4 - E-pts e-lns} }

The e-points and e-lines defined here will be used to construct the
projective extension.
\begin{defn}
\label{Defn. ept. eln.} ~~

\noindent $\bullet$ An\emph{ extension point, }or\emph{ e-point,}
is an equivalence class $\overline{\alpha}$ of pencils of lines.
The e-points $\overline{\alpha}$ and $\overline{\beta}$ are said
to be \emph{equal,} written $\overline{\alpha}=\overline{\beta},$
if $\alpha\approx\beta.$ We say that $\overline{\alpha}$ and $\overline{\beta}$
are \emph{distinct}, written $\overline{\alpha}\ne\overline{\beta},$
if $\alpha'\neq\beta'.$ 

\noindent $\bullet$ An \emph{extension line,} or \emph{e-line,} is
a set $\lambda_{p}$ of e-points, where $p$ is a prime v-line, and
where $\overline{\alpha}\in\lambda_{p}$ if $\overline{\alpha}$ satisfies
these two conditions: 

(1) If $p\neq\varnothing,$ then $p\in\alpha'.$

(2) If $p=\varnothing,$ then $\alpha$ is a parallel pencil.

\noindent We say that the e-lines $\lambda_{p}$ and $\lambda_{q}$
are \emph{equal, }written\emph{ $\lambda_{p}=\lambda_{q},$} or \emph{distinct,}
written $\lambda_{p}\neq\lambda_{q},$ when $p=q,$ or $p\neq q.$
The prime v-line $p$ is called the \emph{root} of $\lambda_{p}.$ 

\noindent $\bullet$ The e-line
\[
\iota=\lambda_{\varnothing}=\{\overline{\alpha}:\alpha\text{ is a parallel pencil}\}
\]
will be called the \emph{line at infinity.}

\noindent $\bullet$ An e-point of the form $\overline{Q^{*}},$ where
$Q\in\mathscr{P},$ will be said to be a\emph{ proper }e-point.\emph{
}An e-point of the form $\overline{\rho},$ where $\rho$ is a regular
pencil, will be called a\emph{ regular }e-point.\emph{ }An e-line
of the form $\lambda_{l},$ where $l\in\mathscr{L},$ will be said
to be a \emph{proper }e-line.\emph{ }\end{defn}
\begin{lem}
\label{Lm-tite-epts-elns} ~

(a) For any e-points $\overline{\alpha}$ and $\overline{\beta},$
if $\neg(\overline{\alpha}\ne\overline{\beta}),$ then $\overline{\alpha}=\overline{\beta}.$ 

(b) For any e-lines $\lambda_{p}$ and $\lambda_{q},$ if $\neg(\lambda_{p}\neq\lambda_{q}),$
then $\lambda_{p}=\lambda_{q}.$ \end{lem}
\begin{proof}
(a) Let $\alpha\subset\rho,$ where $\rho$ is a regular pencil. It
follows from Lemma \ref{Reg vs Equiv} that $\alpha\approx\rho$ and
from Theorem \ref{Thm Union class} that $\alpha'=\rho.$ Let $l\in\beta,$
and suppose that $l\notin\rho;$ then $\alpha'\neq\beta',$ so $\overline{\alpha}\ne\overline{\beta},$
a contradiction. Thus $\neg(l\notin\rho),$ and by Corollary \ref{regular tite}
we have $l\in\rho.$ This shows that $\beta\subset\rho.$ Hence $\alpha\approx\beta,$
and $\overline{\alpha}=\overline{\beta}.$ 

(b) The given condition implies that $\neg(p\neq q),$ so $p\approx q.$
Since these v-lines are prime, it follows that $p=q,$ and hence $\lambda_{p}=\lambda_{q}.$ \end{proof}
\begin{lem}
\label{lambda-l} For any line $l\in\mathscr{L},$ 

(a) $\lambda_{l}=\{\overline{\alpha}:l\in\alpha'\}$ 

(b) $\lambda_{l}\supset\{\overline{Q^{*}}:Q\in l\}\cup\{\overline{l^{*}}\}$ 

\noindent However, equality in (b) would be constructively invalid. \end{lem}
\begin{proof}
Both (a) and (b) are evident. For a Brouwerian counterexample to equality
in (b), let $c\in\mathbb{R},$ let $l$ and $m$ be the lines with
equations $y=0$ and $y=1-cx,$ and construct the pencil $\alpha=\varphi(l,m)$
using Theorem \ref{Thm phi(l,m)}. Since $l\in\alpha,$ we have $\overline{\alpha}\in\lambda_{l}.$
If the hypothesis of equality in (b) is applied to the e-line $\lambda_{l},$
then $\overline{\alpha}$ lies in one of the two indicated sets. If
$\overline{\alpha}=\overline{l^{*}}$, then $\alpha$ is a parallel
pencil, so $c=0.$ If $\overline{\alpha}=\overline{Q^{*}},$ where
$Q=(d,e),$ then $\alpha\subset Q^{*},$ so $e=0$ and $1-cd=0;$
thus $c\ne0.$ Hence LPO results. 
\end{proof}
\noindent The e-point $\overline{l^{*}}$ is called the \emph{tip}
of the e-line $\lambda_{l}$. 
\begin{thm}
An e-line cannot contain fewer than three e-points. \end{thm}
\begin{proof}
Assume that an e-line $\lambda_{p}$ contains fewer than three e-points.
Suppose that $p\neq\varnothing;$ then $p$ is a line, and contains
at least two points. It follows from Lemma \ref{lambda-l} that  $\lambda_{p}$
contains at least two proper e-points, and also the e-point $\overline{p^{*}},$
contradicting the assumption. Thus $p=\varnothing,$ and $\lambda_{p}$
 is the line at infinity $\iota.$ Now $\lambda_{p}$ contains the
e-points corresponding to three mutually distinct parallel pencils
of lines, a final contradiction. \end{proof}
\begin{lem}
If $\overline{\alpha}$ and $\overline{\beta}$ are distinct e-points
on the line at infinity $\iota,$ then $\alpha\cap\beta=\varnothing.$
Conversely, if an e-line $\lambda_{p}$ contains distinct e-points
$\overline{\alpha}$ and $\overline{\beta},$ where $\alpha$ and
$\beta$ are complete pencils with $\alpha\cap\beta=\varnothing,$
then $\lambda_{p}=\iota.$ \end{lem}
\begin{proof}
The first statement follows from Theorem \ref{Thm Par Pens }. For
the second statement, first note that it follows from Corollary \ref{Cor parpen}
that $\alpha$ and $\beta$ are parallel pencils, and then from Theorem
\ref{Thm Pens } that $\alpha'$ and $\beta'$ are also parallel pencils.
Now suppose that $p\neq\varnothing;$ then $p\in\alpha'\cap\beta',$
and it follows from Theorem \ref{Thm Par Pens } that $\alpha'\approx\beta',$
so $\overline{\alpha}=\overline{\beta},$ a contradiction. Hence $p=\varnothing,$
and $\lambda_{p}=\iota.$ 
\end{proof}
\noindent \emph{Note on the stem of an e-line. }As noted at the beginning
of Section \ref{> Section 3 - V-lines}, the notion of \emph{v-line}
arose in connection with the set of proper e-points that lie on an
e-line. For any e-line $\lambda_{p},$ define
\[
s_{p}=\{Q\in\mathscr{P}:\overline{Q^{*}}\in\lambda_{p}\}
\]
 and call this the \emph{stem} of $\lambda_{p}.$ 

It is easily seen that $p\subset s_{p},$ that $\neg(p\neq s_{p}),$
and that $p=s_{p}$ when $p\ne\varnothing.$ However, the statement
\emph{{}``The stem of any e-line is a v-line''} is constructively
invalid. For a Brouwerian counterexample, let $c\in\mathbb{R},$ and
set $q=\{(t,0):t\in\mathbb{R}\text{ and }c=0\}\cup\{(0,t):t\in\mathbb{R}\text{ and }c\ne0\}.$
It is clear that $q$ is a v-line; construct the e-line $\lambda_{p}$
with root $p=q'.$ To show that $\overline{O^{*}}\in\lambda_{p},$
where $O$ is the origin, note first that the condition $p=\varnothing$
is ruled out by Lemma \ref{v-line class}, so we need consider only
the condition $p\ne\varnothing.$ Suppose that $O\notin p,$ and suppose
further that $c\neq0.$ Then $p=q=m_{0},$ where $m_{0}$ is the $y$-axis,
so $O\in p,$ a contradiction. Thus $c=0,$ and $p=q=l_{0},$ a contradiction.
Thus $O\in p,$ and $p\in O^{*}.$ Since $O^{*}$ is a prime pencil,
this shows that $\overline{O^{*}}\in\lambda_{p};$ thus $O\in s_{p},$
so $s_{p}$ is nonvoid. Note that if $c=0,$ then $s_{p}=p=q=l_{0},$
while if $c\neq0;$ then $s_{p}=p=q=m_{0}.$ By hypothesis, $s_{p}$
is a line; thus either $s_{p}\neq l_{0},$ or $s_{p}\neq m_{0}.$
It follows that either $\neg(c=0),$ or $c=0.$ Hence WLPO results.

\section{Projective extension\label{> Section 5 - Proj ext} }

The extension will be constructed and the main projective axioms will
be verified. 
\begin{defn}
~

\noindent $\bullet$ We denote by $\mathscr{P^{*}}$ the family of
all e-points, and by $\mathscr{L^{*}}$ the family of all e-lines,
retaining the equality and inequality relations, and the relation
$\overline{\alpha}\in\lambda_{p},$ adopted in Definition \ref{Defn. ept. eln.}. 

\noindent $\bullet$ For any e-point $\overline{\alpha},$ and any
e-line $\lambda_{p}$, we say that $\overline{\alpha}$ \emph{lies}
\emph{outside} $\lambda_{p}$, written $\overline{\alpha}\notin\lambda_{p},$
if $\overline{\alpha}\neq\overline{\rho}$ for all regular e-points
$\overline{\rho}$ in $\lambda_{p}.$ 

\noindent $\bullet$  $\mathscr{G^{*}}=(\mathscr{P^{*}},\mathscr{L^{*}})$
is the \emph{projective extension} of the incidence plane $\mathscr{G=}(\mathscr{P,L})$. \end{defn}
\begin{thm}
Let $\overline{\alpha}$ be any e-point, and $\lambda_{p}$ any e-line.
If $\neg(\overline{\alpha}\notin\lambda_{p}),$ then $\overline{\alpha}\in\lambda_{p}.$ \end{thm}
\begin{proof}
First let $p\neq\varnothing,$ and suppose that $p\notin\alpha'.$
If $\overline{\rho}$ is any regular e-point on $\lambda_{p},$ then
$p\in\rho,$ so $\alpha'\neq\rho,$ and $\overline{\alpha}\neq\overline{\rho};$
this shows that $\overline{\alpha}\notin\lambda_{p},$ a contradiction.
Thus we have $\neg(p\notin\alpha'),$ and it follows from Theorem
\ref{Thm. tite} that $p\in\alpha'.$ 

Now let $p=\varnothing;$ thus $\lambda_{p}=\iota,$ the line at infinity.
Let $l,m\in\alpha,$ and suppose that $l$ intersects $m.$ Let $\overline{\rho}$
be any regular e-point on $\lambda_{p};$ then $\rho=n^{*}$ for some
line $n.$ By Axiom L2 of {[}M07{]} we may assume that $n$ intersects
$l,$ so $l\notin n^{*}.$ Thus $\alpha'\neq n^{*},$ and $\overline{\alpha}\neq\overline{\rho}.$
This shows that $\overline{\alpha}\notin\lambda_{p},$ a contradiction.
Thus $l\parallel m,$ and this shows that $\alpha$ is a parallel
pencil. 

Hence $\overline{\alpha}\in\lambda_{p}.$\end{proof}
\begin{thm}
\label{com line} There exists a unique e-line passing through any
two distinct e-points. \end{thm}
\begin{proof}
(i) Let $\overline{\alpha}$ and $\overline{\beta}$ be distinct e-points;
we may assume that $\alpha$ and $\beta$ are complete pencils. By
Lemma \ref{Lm. core is quasi-line}, the core $r=\alpha\sqcap\beta$
is a v-line; set $p=r'.$ First let $p\ne\varnothing,$ and assume
that $p\notin\alpha'.$ Suppose that $r\neq\varnothing;$ then $p=r\in\alpha,$
a contradiction. Thus $r=\varnothing,$ so also $p=\varnothing,$
a contradiction. This shows that $\neg(p\notin\alpha'),$ and it follows
from Theorem \ref{Thm. tite} that $p\in\alpha'.$ Now let $p=\varnothing;$
then also $r=\varnothing.$ Since $\alpha$ and $\beta$ are complete
pencils, it follows from Corollary \ref{Cor parpen} that $\alpha$
is a parallel pencil. Thus $\overline{\alpha}\in\lambda_{p},$ and
similarly $\overline{\beta}\in\lambda_{p}.$ 

(ii) Let $\lambda_{p}$ and $\lambda_{q}$ be e-lines, each passing
through two distinct e-points $\overline{\alpha}$ and $\overline{\beta},$
and assume that $p\neq q.$ Suppose that $p\neq\varnothing;$ then
$\alpha'\cap\beta'=\{p\}.$ Suppose further that $q\neq\varnothing;$
then also $\alpha'\cap\beta'=\{q\},$ so $p=q,$ a contradiction.
Thus $q=\varnothing.$ Now $\alpha$ and $\beta$ are parallel pencils,
and so also are $\alpha'$ and $\beta'.$ It follows from Theorem
\ref{Thm Par Pens } that $\alpha'\cap\beta'=\varnothing;$ a contradiction.
Thus $p=\varnothing.$ Similarly, $q=\varnothing,$ so $p=q,$ a final
contradiction. This shows that $\neg(p\neq q),$ so $p\approx q.$
Since $p$ and $q$ are prime v-lines, it follows that $p=q,$ and
hence $\lambda_{p}=\lambda_{q}.$ \end{proof}
\begin{cor}
\label{Cor Inv Core} \emph{Invariance of the core.} Let $\lambda_{p}$
be an e-line, and let \textup{\emph{$\overline{\alpha}$ and $\overline{\beta}$}}
be any distinct e-points on $\lambda_{p},$ where \textup{$\alpha$
and $\beta$} are prime pencils.\textup{\emph{ Then }}$\alpha\sqcap\beta=p.$ \end{cor}
\begin{proof}
Set $r=\alpha\sqcap\beta,$ and $q=r'.$ With the construction of
Theorem \ref{com line}, $\overline{\alpha}$ and $\overline{\beta}$
will lie on the e-line $\lambda_{q}.$ Thus $p=q;$ this shows that
$r\subset p.$ If $Q\in p,$ then $p\neq\varnothing,$ so $p\in\alpha\cap\beta.$
It follows from Theorem \ref{Thm Pens } that $p=r,$ and thus $Q\in r;$
this shows that $p\subset r.$ 
\end{proof}
\noindent \emph{Note for Corollary \ref{Cor Inv Core}.} A stronger
statement, without the condition {}``prime'' pencils, would be constructively
invalid. For a Brouwerian counterexample, let $c$ be a real number
with $\neg(c=0),$ set $\alpha=(0,0)^{*},$ and set $\beta=\{l_{t}:t\in\mathbb{R}\text{ and }c\ne0\}.$
Then $\overline{\alpha},\overline{\beta}\in\lambda_{l_{0}},$ and
$\alpha\sqcap\beta=\{(t,0):t\in\mathbb{R}\text{ and }c\neq0\}.$ By
hypothesis, we have $\alpha\sqcap\beta=l_{0},$ so $c\neq0.$ Thus
MP results. 

\bigskip{}

\noindent Classically, the projective axiom concerning a common point
for any two lines need be verified only for distinct lines, there
being no reason to consider identical lines. Constructively, however,
there are always innumerable pairs of lines for which we do not know,
at present, whether they are identical or distinct. Thus we require
a theorem that deals with two arbitrary lines. The point obtained
will be common to the lines in any eventuality, allowing for any possible
future discovery that the lines are distinct, or are identical. 
\begin{thm}
Any two e-lines pass through a common e-point. If the e-lines are
distinct, then the common e-point is unique. \end{thm}
\begin{proof}
Given the e-lines $\lambda_{p}$ and $\lambda_{q},$ construct the
complete pencil $\gamma=\varphi(p,q),$ using Theorem \ref{Thm. phi(p,q)}.
If $p\neq\varnothing,$ then $p\in\varphi_{0},$ so $p\in\gamma\subset\gamma'.$
If $p=\varnothing,$ then it follows from Lemma \ref{Lemma for phi}
that $\gamma$ is a parallel pencil. Hence $\overline{\gamma}\in\lambda_{p},$
and similarly $\overline{\gamma}\in\lambda_{q}.$ Thus $\overline{\gamma}$
is a common e-point. If $\lambda_{p}\neq\lambda_{q},$ then uniqueness
of the common e-point follows from Theorem \ref{com line} and Lemma
\ref{Lm-tite-epts-elns}.
\end{proof}
\noindent \emph{Note on cotransitivity.} The Brouwerian counterexample
in {[}M11{]} shows that in any projective-type extension, cotransitivity
for points is constructively incompatible with the existence of a
common point for any two lines. It follows that the statement {}``In
the extension $\mathscr{G^{*}},$ the cotransitivity property holds
for e-points'' is constructively invalid. 

For a direct Brouwerian counterexample, using the constructions of
$\mathscr{G^{*}},$ let $c\in\mathbb{R},$ set $d=\max\{c,0\},$ set
$e=\max\{-c,0\},$ let $l$ be the line on $\mathbb{R}^{2}$ with
equation $y=1-dx,$ let $m$ be the line with equation $x=1-ey,$
set $\alpha=\varphi(l_{0},l),$ set $\beta=\varphi(m_{0},m),$ where
$m_{0}$ is the $y$-axis, set $p=\{(1/c,t):t\in\mathbb{R}\,\,\text{and}\,\, c>0\},$
set $q=\{(t,1/c):t\in\mathbb{R}\,\,\text{and}\,\, c<0\},$ and set
$\gamma=\varphi(p,q)$ using Theorem \ref{Thm. phi(p,q)}. To show
that $\overline{\alpha}\neq\overline{\beta},$ it will suffice to
show that $l_{0}\notin\beta'.$ For any line $n\in\beta',$ it follows
from Axiom L2 in {[}M07{]} that either $n$ intersects $l_{0}$ or
$n$ intersects $m_{0}.$ In the first case, $n\neq l_{0}.$ In the
second case, say $n\times m_{0}=(0,h);$ then $\beta'=(0,h)^{*}.$
Since $m\in\beta,$ we have $0=1-eh,$ so $h\neq0,$ and thus $(0,h)\notin l_{0};$
it follows that $n\neq l_{0}.$ This shows that $l_{0}\notin\beta'.$
Thus $\alpha'\neq\beta',$ and $\overline{\alpha}\neq\overline{\beta}.$
By hypothesis, we have either $\overline{\gamma}\neq\overline{\alpha}$
or $\overline{\gamma}\neq\overline{\beta}.$ In the first case, suppose
that $c<0.$ Then $d=0,$ so $\alpha=l_{0}^{*}.$ Also, $p=\varnothing,$
and $q=l_{1/c},$ so $\gamma=q^{*}=l_{0}^{*}=\alpha,$ a contradiction;
thus $c\geq0.$ Similarly, in the case $\gamma\neq\overline{\beta}$
we find that $c\leq0.$ Hence LLPO results. 

For an alternative counterexample, use the v-line $p=\{(t,0):t\in\mathbb{R}\,\,\text{and}\,\, c=0\}\cup\{(0,t):t\in\mathbb{R}\,\,\text{and}\,\, c\ne0\}$
to construct the pencil $\gamma=\varphi(p,p).$ Assuming cotransitivity
for e-points, we have either $\overline{\gamma}\neq\overline{l_{0}^{*}}$
or $\overline{\gamma}\neq\overline{m_{0}^{*}},$ where $m_{0}$ is
the $y$-axis. Hence WLPO results. 

Cotransitivity for e-lines in the extension $\mathscr{G^{*}}$ is
also constructively invalid; this may be seen using the counterexample
at the end of Section \ref{> Section 4 - E-pts e-lns}.
\begin{thm}
The projective plane $\mathscr{G^{*}}=(\mathscr{P^{*}},\mathscr{L^{*}})$
is an extension of the incidence plane $\mathscr{G=}(\mathscr{P,L}).$ \end{thm}
\begin{proof}
Set $\mathscr{P'}=\{\overline{Q^{*}}:Q\in\mathscr{P}\},$ the family
of proper e-points, and set $\mathscr{L'}=\{\lambda_{l}:l\in\mathscr{L}\},$
the family of proper e-lines. Then the mappings $Q\mapsto\overline{Q^{*}}$
$(Q\in\mathscr{P),}$ and $l\mapsto\lambda_{l}$ $(l\in\mathscr{L}),$
map $\mathscr{G}$ into $\mathscr{G^{*},}$ with image $\mathscr{G^{'}}=(\mathscr{P^{'}},\mathscr{L^{'}}).$ 

If $Q=R,$ then clearly $\overline{Q^{*}}=\overline{R^{*}}.$ Conversely,
if $\overline{Q^{*}}=\overline{R^{*}},$ then $Q^{*}\approx R^{*},$
and since these pencils are regular, we have $Q^{*}=R^{*}.$ Select
distinct lines $l$ and $m$ in $Q^{*};$ thus $l\times m=Q$. Since
also $l,m\in R^{*},$ we have $l\times m=R,$ so $Q=R.$ 

If $Q\neq R,$ we may select a line $l$ so that $Q\in l$ but $R\notin l.$
Thus $l\in Q^{*},$ but $l\notin R^{*},$ so $Q^{*}\neq R^{*};$ since
these pencils are prime, we have $\overline{Q^{*}}\neq\overline{R^{*}}.$
Conversely, if $\overline{Q^{*}}\neq\overline{R^{*}},$ then $Q^{*}\neq R^{*},$
and we may select a line $m\in Q^{*}$ with $m\notin R^{*}.$ Then
$Q\in m$ but $R\notin m,$ so $Q\neq R.$ 

For any lines $l$ and $m,$ it follows from Definition \ref{Defn. ept. eln.}
that $l=m$ (or $l\neq m)$ if and only if $\lambda_{l}=\lambda_{m}$
(or $\lambda_{l}\neq\lambda_{m}).$ Thus the mappings preserve the
equality and inequality relations for points and lines. 

Clearly, $Q\in l$ if and only if $\overline{Q^{*}}\in\lambda_{l}.$
Now let $Q\notin l.$ For any e-point $\overline{\beta}\in\lambda_{l},$
we have $l\in\beta',$ but $l\notin Q^{*};$ thus $Q^{*}\neq\beta',$
so $\overline{Q^{*}}\neq\overline{\beta}.$ Thus $\overline{Q^{*}}\notin\lambda_{l}.$
Conversely, let $\overline{Q^{*}}\notin\lambda_{l}.$ For any point
$R\in l,$ we have $\overline{R^{*}}\in\lambda_{l};$ thus $\overline{Q^{*}}\neq\overline{R^{*}},$
so $Q\neq R.$ Thus $Q\notin l$. Thus the mappings preserve the relations
\emph{point on a line }and \emph{point outside a line.} 

Hence the mappings $Q\mapsto\overline{Q^{*}}$ and $l\mapsto\lambda_{l}$
form a strict isomorphism of $\mathscr{G}$ onto $\mathscr{G}',$
a sub-plane of $\mathscr{G^{*}}.$ 
\end{proof}

\section*{References\label{Section - References} }

\noindent {[}B65{]} E. Bishop, Book Review: \emph{The Foundations
of Intuitionistic Mathematics, }by S. C. Kleene and R. E. Vesley.
Bull. Amer. Math. Soc. 71 (1965), no. 6, 850--852.  

\noindent {[}B67{]} ------, \emph{Foundations of Constructive Analysis.}
McGraw-Hill Book Co., New York-Toronto-London, 1967.

\noindent {[}BB85{]} E. Bishop and D. Bridges, \emph{Constructive
Analysis.} Springer-Verlag, Berlin, 1985.

\noindent {[}BR87{]} D. Bridges and F. Richman, \emph{Varieties of
Constructive Mathematics.} Cambridge University Press, Cambridge,
1987.

\noindent {[}Br08{]} L. E. J. Brouwer, De onbetrouwbaarheid der logische
principes. Tijdschrift voor Wijsbegeerte 2 (1908), 152-158. English
translation, \textquotedblleft{}The Unreliability of the Logical Principles\textquotedblright{},
pp. 107\textendash{}111, in Heyting, A. (ed.), \emph{L. E. J. Brouwer:
Collected Works 1: Philosophy and Foundations of Mathematics.} Amsterdam-New
York, Elsevier, 1975.

\noindent {[}D63{]} D. van Dalen, Extension problems in intuitionistic
plane projective geometry I, II, Indag. Math. 25 (1963), 349-383. 

\noindent {[}H59{]} A. Heyting, Axioms for intuitionistic plane affine
geometry, in L. Henkin, P. Suppes, A. Tarski (eds), \emph{The Axiomatic
Method, with special reference to geometry and physics: Proceedings
of an international symposium held at the University of California,
Berkeley, December 26, 1957 - January 4, 1958.} North-Holland, Amsterdam,
1959, 160-173. 

\noindent {[}M83{]} M. Mandelkern, \emph{Constructive Continuity.}
Memoirs Amer. Math. Soc. 42 (1983), nr. 277. 

\noindent {[}M85{]} ------, Constructive mathematics, Math. Mag. 58
(1985), 272-280. 

\noindent {[}M88{]} ------, Limited omniscience and the Bolzano-Weierstrass
principle, Bull. London Math. Soc. 20 (1988), 319-320. 

\noindent {[}M89{]} ------, Brouwerian counterexamples, Math. Mag.
62 (1989), 3-27. 

\noindent {[}M07{]} ------, Constructive coordinatization of Desarguesian
planes, Beitr�ge Algebra Geom. 48 (2007), 547-589.

\noindent {[}M11{]} ------, The common point problem in constructive
projective geometry, preprint, \url{http://www.zianet.com/mandelkern/}. 

\noindent {[}R82{]} F. Richman, Meaning and information in constructive
mathematics, Amer. Math. Monthly 89 (1982), 385-388. 

\noindent {[}R02{]} ------, Omniscience principles and functions of
bounded variation, Math. Logic Quart. 42 (2002), 111-116. 

\noindent {[}S70{]} G. Stolzenberg, Review of E. Bishop, \emph{Foundations
of Constructive Analysis,} Bull. Amer. Math. Soc. 76 (1970), 301-323.
\\

\noindent New Mexico State University

\noindent Las Cruces, New Mexico, USA

\noindent 2012 February 2 
\end{document}